\documentclass{amsart}

\usepackage{latexsym}
\usepackage{amssymb}

\newtheorem{THM}{Theorem}
\newtheorem{Lemma}[THM]{Lemma}

\renewcommand{\Re}{\mathbb{R}}
\newcommand{\be}{\begin{equation}}
\newcommand{\ee}{\end{equation}}

\newcommand{\ra}{\rightarrow}

\newcommand{\lp}{\left(}
\newcommand{\rp}{\right)}
\newcommand{\lb}{\left[}
\newcommand{\rb}{\right]}
\newcommand{\lc}{\left\{}
\newcommand{\rc}{\right\}}
\newcommand{\lab}{\left|}
\newcommand{\rab}{\right|}
\newcommand{\Lap}{\Delta}

\newcommand{\sF}{\mathcal{F}}

\newcommand{\E}{\mathbb{E}}
\newcommand{\Prob}{\mathbb{P}}

\DeclareMathOperator{\Area}{Area}

\renewcommand{\phi}{\varphi}

\begin{document}

\title{On parabolicity and area growth of minimal surfaces}
\author{Robert W.\ Neel}
\address{Department of Mathematics, Lehigh University, Bethlehem, PA}
\begin{abstract}We establish parabolicity and quadratic area growth for minimal surfaces-with-boundary contained in regions of $\Re^3$ which are within a sub-logarithmic factor of the exterior of a cone.  Unlike previous work showing that these two properties hold for minimal surfaces-with-boundary contained between two catenoids, we do not make use of universal superharmonic functions.  Instead, we use stochastic methods, which have the additional feature of giving a type of parabolicity in a more general context than Brownian motion on a minimal surface.
\end{abstract}
\email{robert.neel@lehigh.edu}
\date{February 15, 2011}

\maketitle

\section{Introduction}

We are interested in controlling the geometry and conformal structure of minimal surfaces-with-boundary in $\Re^3$, under the assumption that they are contained in certain rotationally symmetric regions.  In particular, let $(x_1, x_2, x_3)$ be standard Euclidean coordinates on $\Re^3$ and let $r=\sqrt{x_1^2+x_2^2}$.   We consider minimal surfaces-with-boundary contained in sets of the form $\{ |x_3| \leq f(r)\text{ and }r>C\}$ for some positive, increasing, continuous function $f$ and some positive constant $C$.  We are interested in two possible properties of such surfaces, parabolicity and quadratic area growth.  A surface-with-boundary is parabolic if the boundary is non-empty and any bounded harmonic function on the surface is determined by its boundary values.  (This is not the standard notion of parabolicity for a Riemann surface, but it will be a convenient usage for us.)  Equivalently, a surface is parabolic if Brownian motion on the surface, started from any interior point, hits the boundary (in finite time) almost surely.  Now for positive $\rho$, let $B_{\rho}=\{x_1^2+x_2^2+x_3^2 \leq \rho^2\}$ be the ball of radius $\rho$.  We say that a surface $M$ has quadratic area growth if $M\cap B_{\rho}$ has area less than or equal to $C \rho^2$ for sufficiently large $\rho$, for some positive constant $C$.  (Note that this is an extrinsic notion of quadratic area growth.)

We prove that both of these properties hold for minimal surfaces-with-boundary contained in regions that are within sub-logarithmic factors of the exterior of a cone (that is, $f(r)$ as above is within a sub-logarithmic factor of being linear), under mild additional assumptions like stochastic completeness or properness of the immersion.  These results are given in Theorems \ref{THM:Para} and \ref{THM:Area}.  We prove both results using a stochastic approach, in contrast to the analytic approach via universal superharmonic functions (as described below).  One additional feature of this approach is that parabolicity can be proven for a broader class of processes than Brownian motion on minimal surfaces, in a sense described in Theorem \ref{THM:MartPara}.

In the remainder of this section, we give additional background on this problem and make a few preliminary observations.  We say that a surface $M$ is contained between two catenoids if $M$ is contained in the region $\{|x_3|\leq c \log r\text{ and }r\geq 2\}$
for some $c>0$.  In \cite{CKMR}, parabolicity and quadratic area growth were shown for properly immersed minimal surfaces with compact boundary contained between two catenoids (see also Section 6 of the comprehensive survey \cite{MPSurvey}).  This was done using the technique of universal superharmonic functions, where a universal superharmonic function is a function on a subset of $\Re^3$ with the property that its restriction to any minimal surface is superharmonic.  More specifically, the method of \cite{CKMR} is to find a universal superharmonic function on the region between two catenoids which is also proper and bounded from below.  The existence of such a function immediately shows that any minimal surface of this type is parabolic.  Further, using that the gradient and Laplacian of their universal superharmonic function are natural geometric quantities, they were able to deduce quadratic area growth.  Bill Meeks has asked whether universal superharmonic functions which are proper and bounded from below can be found on the exterior of sufficiently large cones, mainly for the purpose of extending the parabolicity and quadratic area growth results just mentioned.  The results of the present paper bring us much closer (within sub-logarithmic factors) to the region conjectured by Meeks (of course, they leave open the question of the existence of universal superharmonic functions themselves on such regions).

We also mention the complimentary project of finding regions (preferably as small as possible) where one can minimally (conformally) immerse a surface of arbitrary conformal structure.  Recent progress in this direction can be found in \cite{AlarconLopez}.

The motivation for studying minimal surfaces-with-boundary restricted to this type of region comes in part from the fact that such surfaces arise as (representatives of) ends of complete minimal surfaces.  In particular, the results of \cite{CKMR} on minimal surfaces-with-boundary contained between two catenoids, mentioned above, were used in the same paper in order to advance the classification of properly embedded minimal surfaces by controlling the behavior of certain ends.  In light of this connection, as well as the fact that parabolicity and quadratic area growth are global properties, we are primarily interested in the large $r$ behavior of the surfaces, and thus our results are formulated without much concern for the small $r$ behavior of the surfaces.  In this context, we note that parabolicity and quadratic area growth are preserved under homothetic rescaling (which also preserves minimality of the surface) and under removing compact subsets of the surface.  This allows the results proven below to be applied somewhat more broadly than the hypotheses might otherwise indicate.

We make one more preliminary observation.  Suppose that $M$ is contained in the exterior of a cone, that is, $M$ is contained in the region $\{(x_1,x_2,x_3): |x_3| \leq r \}$ for some positive constant $c$.  Let $T_{\rho}=\{(x_1,x_2,x_3): r \leq \rho\}$ for positive $\rho$ (so that $T_{\rho}$ is a cylindrical region).  Then we see that $M$ has quadratic area growth if and only if the area of $M\cap T_{\rho}$ is bounded by $C\rho^2$ for large $\rho$, for some positive constant $C$.  That is to say, quadratic area growth with respect to the distance from the origin is equivalent to quadratic area growth with respect to $r$ for such surfaces.  All of the surfaces we consider in this paper are contained in the exterior of a cone and thus satisfy this equivalence.  In what follows, we will frequently use the condition on $M\cap T_{\rho}$ as our criterion for quadratic area growth without further comment.

The author would like to thank Bill Meeks for introducing him to this problem in the context of Oberwolfach's Arbeitsgemeinschaft on minimal surfaces held in October 2009.  In particular, this paper addresses several questions that the author raised in his talk there.  The author would also like to thank Bennett Eisenberg for advice about Markov chains, and Graham Smith and Rob Kusner for discussions about universal superharmonic functions and minimal ends.

\section{Background on Markov chains}\label{Sect:MC} 

A central tool for us will be the comparison of Brownian motion on $M$ (or a more general martingale) to a Markov chain.  In preparation for this, we briefly review some relevant facts about Markov chains, formulated for our case of interest.

We take a subset of the non-negative integers, $\Omega_L=\{L, L+1,L+2, \ldots\}$ for some non-negative integer $L$, as our state space.  We are interested in discrete time Markov chains $(Y_n, n\in\{0, 1,\ldots\})$ that evolve in the following way.  If $Y_n=m>L$, then $Y_{n+1}$ is $m+1$ with probability $p_m$ and $m-1$ with probability $q_m=1-p_m$.  We stop the process upon hitting the boundary at $L$.  So $Y_n$ is a time-homogenous nearest-neighbor random walk on $\Omega_L$, stopped at $L$.  We will always assume that $p_m>0$ and $q_m>0$ for all $m>L$, so that our chain is irreducible.  Obviously, one could reduce this to the usual case of $L=0$ by a simple translation.  However, this more general set-up will make the notation somewhat less burdensome later on.

For $m\geq L+1$, we let
\[
A_m= 1+\frac{q_{L+1}}{p_{L+1}} +\frac{q_{L+1}q_{L+2}}{p_{L+1}p_{L+2}} +\cdots +\frac{q_{L+1}\cdots q_m}{p_{L+1}\cdots p_m} .
\]
A well-known computation (see, for example, Theorems 5.3.7 and 5.3.8 of \cite{DurrettProb}) shows that $A_m$ is the reciprocal of the probability that the process, starting from $L+1$, hits $m+1$ before $L$ (we note that, under our assumptions, the process almost surely hits either $L$ or $m+1$ in finite time).  Thus, the chain is parabolic, meaning that it almost surely hits $L$, if and only if $A_m\ra\infty$ as $m\ra\infty$.

To give an example which will be relevant later, let $q_m=m/(2m+1)$ and $p_m=(m+1)/(2m+1)$.  Then 
\[
\frac{q_{L+1}\cdots q_m}{p_{L+1}\cdots p_m} = \frac{1}{m} \quad\text{and}\quad A_m = \sum_{k=L}^m\frac{1}{k} .
\]
Because the harmonic series is divergent, it follows that this chain is parabolic, and more generally, so is any chain with $p_m\leq (m+1)/(2m+1)$ for large $m$.  Further, we think of this chain as a borderline case, in the sense that the harmonic series is a borderline case for the divergence of series.

In addition, we will need a result on an expected number of upcrossings.  We assume that the chain is parabolic.  Then if we start the process at $Y_0=m$, we let $\sigma_1$ be the first hitting time of $m+1$, which may be infinite if the process reaches $L$ first.  Then we let $\tilde{\sigma}_1$ be the first hitting time of $m$ after $\sigma_1$, $\sigma_2$ the first hitting time of $m+1$ after $\tilde{\sigma}_1$, and so on.  Note that if $\sigma_i$ is finite, then $\tilde{\sigma}_i$ is also finite, almost surely, because the chain is parabolic.  If we let $u_m = \sup\{ m: \sigma_m<\infty \}$ (with the convention that $u_m=0$ if $\sigma_1=\infty$), then $u_m$ is the number of upcrossings from $m$ to $m+1$.  Because the chain is parabolic, $u_m$ is almost surely finite.  A straightforward adaptation of the computation that leads to $A_m$ shows that
\[\begin{split}
& \Prob\lp \text{$Y_n$ hits $L$ before $m+1$} |Y_0=m\rp \\
&\quad\quad\quad= \lp 1+\frac{p_m}{q_m}+\frac{p_m p_{m-1}}{q_m q_{m-1}}+\cdots+\frac{p_m\cdots p_{L+1}}{q_m\cdots q_{L+1}}\rp^{-1} .
\end{split}\]
Further, by the Markov property, the number of upcrossings from $m$ to $m+1$ is a geometric random variable with the above as its parameter, and thus its expectation is
\[
\E\lp u_m| Y_0=m \rp = 1+\frac{p_m}{q_m}+\frac{p_m p_{m-1}}{q_m q_{m-1}}+\cdots+\frac{p_m\cdots p_{L+1}}{q_m\cdots q_{L+1}} .
\]

\section{Control of the radial process}\label{Sect:Radial}

Having recalled some basic facts about Markov chains, we now discuss how they arise in the present work.    We begin by introducing a broader class of processes, which we will use in our discussion of parabolicity.  Consider an $\Re^3$-martingale $Z_t$ which solves an (Ito) SDE of the form
\[
dZ_t = \sigma_t \, dW^{\Re^3}_t , \quad \text{where $W^{\Re^3}_t$ is an $\Re^3$-Brownian motion,}
\]
and where $\sigma_t$ is an adapted process taking values in the set of $3\times3$ matrices conjugate under rotation to
\[
\begin{bmatrix}
1 & 0 & 0 \\
0 & 1 & 0 \\
0 & 0 & 0
\end{bmatrix} .
\]
That is, $\sigma_t$ takes values in the set of projections onto planes.  For simplicity, we will also assume that $\sigma_t$ is a continuous process.  Intuitively, such a $Z_t$ is infinitesimally Brownian motion on the plane corresponding to $\sigma$, where $\sigma$ is allowed to vary.  One possibility is for $\sigma_t$ to be projection onto the tangent plane of a minimal surface.  More specifically, let $M$ be a stochastically complete minimal surface-with-boundary.  Stochastic completeness means that Brownian motion on $M$ exists until it hits the boundary (that is, it doesn't explode), and we note that all properly immersed minimal surfaces-with-boundary are stochastically complete.  Then $Z_t$ is Brownian motion on $M$, immersed into $\Re^3$.  (Here we are using the fact that such a $Z_t$ stays on the minimal surface, so $\sigma_t$ is well-defined as projection onto the tangent plane at $Z_t$, where self-intersections are handled in the usual way.)  This is the motivating case.  However, the point is that nothing in our approach to parabolicity will require the plane field corresponding to $\sigma_t$ to be integrable, so that $Z_t$ need not be restricted to a surface.  In particular, in what follows it is appealing to consider the case when we choose $\sigma_t$ as a control in order to try to force $Z_t$ to remain in some subset of $\Re^3$.  As usual, we also wish to allow the possibility that our process $Z_t$ is stopped in finite time; in particular, $Z_t$ runs until some stopping time $\zeta$, which may be infinite.  In the minimal surface-with-boundary case, we will always take $\zeta$ to be the hitting time of the boundary.

For convenience, we call such a $Z_t$ a rank 2 martingale.  The terminology is motivated by the fact that $\sigma_t$ has rank 2, which means that the process infinitesimally evolves along a plane.  Of course, we're assuming more than just that $\sigma_t$ has rank 2, corresponding to assuming that $X_t$ is not merely infinitesimally evolving along a plane but is doing so ``like Brownian motion.''  Nonetheless, making our terminology more specific at the expense of making it more cumbersome doesn't seem worth it.

Recall that we let $r=\sqrt{x_1^2+x_2^2}$ be the radial coordinate, in the sense of cylindrical coordinates.  We will assume that $Z_t$ is contained in the region $\{r\geq e^L\}$ for some non-negative integer $L$, and that it is stopped if it ever hits $\{r= e^L\}$.  For convenience, we will also assume that $r(Z_0)=e^{L+1}$.  We now discretize $Z_t$ (in particular, its radial component $r_t=r(Z_t)$) in the following way.  Set $\sigma_0=0$, and let $\sigma_1$ be the first time $\log r_t$ hits either $L$ (in which case the process is stopped) or $L+2$.  Then we inductively define $\sigma_{n+1}$ for $n\geq 1$ to be the first time $\log r_t$ hits either $\log r_{\sigma_{n}}+1$ or $\log r_{\sigma_{n}}-1$.  As usual, we allow $\sigma_n$ to be $\infty$ if the event never occurs.

Let $n_t$ be a unit vector in the kernel of $\sigma_t$, and let $n_{3,t}$ be its $x_3$-component.  The notation is chosen because when $Z_t$ is Brownian motion on a minimal surface, $n_t$ is a unit normal to $M$ (at the current point).  Then a simple application of Ito's rule shows that $\alpha_t=1-n_{3,t}^2$ is the infinitesimal drift of $r^2_t$, meaning that $r(Z_t)^2 - \int_0^t \alpha_t \, dt$ is a martingale.  In the minimal surface case, this is just the fact that $(\Lap_M/2) r^2 = 1-n_3^2$, where $n_3$ is the $x_3$-component of the unit normal to $M$.  Since $1\leq \alpha_t\leq 2$, it follows that the expectation of $r^2_t$ grows at a rate between 1 and 2.  In particular, we see that $\log r_t$ almost surely leaves any interval of the form $(m+1,m-1)$ for $m\geq L+1$, unless the process is stopped first.  Thus, almost surely, $\sigma_n$ will be $\infty$ if and only $\zeta$ is finite and occurs before $\sigma_n$.  So we have that $X_n=\log r_{\sigma_n}$ is a nearest-neighbor random walk on $\{L,L+1,\ldots \}$, started at $X_0=L+1$ and stopped at the random time $\zeta$ (where this has the obvious meaning in terms of the underlying process $Z_t$).  If the process is stopped before $n$, in the sense that $\zeta<\sigma_n$, we will say that $X_n$ does not exist.  Note that $X_n$ is not, in general, Markov.  This is true even in the case when $Z_t$ is Brownian motion on a minimal surface-with-boundary $M$, since the transition probabilities (as well as $\zeta$) depend on where $Z_{\sigma_n}$ is (on $M$).

We now further assume that $Z_t$ is contained in the region
\[
\{|x_3|\leq f(r) \}
\]
for some continuous, non-negative, monotone non-decreasing function $f:[e^L,\infty)\ra[0,\infty)$.  Our goal is to find a Markov chain with which to compare $X_n$.  The first step is to find upper bounds on $\Prob\lp X_{n+1}=m+1 | X_n=m\rp$ that depend only on $\{X_n=m\}$ and not on anything else in $\sF^Z_{\sigma_n}$ (here $\sF^Z_{\sigma_n}$ is the $\sigma$-algebra generated by $Z_t$ up to time $\sigma_n$).

Recall the following basic estimate for minimal surfaces (see Lemma 2.3 of \cite{CKMR})
\begin{Lemma}\label{Lem:Basic}
Let $M$ be a minimal surface-with-boundary and assume that $r\neq 0$ on $M$.  Then $|\Lap_M \log r| \leq |\nabla_M x_3|^2/r^2$ (on the interior of $M$).
\end{Lemma}
The key observation for us is that this generalizes immediately to rank 2 martingales.  We let $\beta_t$ and $\gamma_t$ be the infinitesimal drifts of $\log r_t$ and $x_{3,t}^2 = x^2_3(Z_t)$, respectively, analogously to our earlier introduction of $\alpha_t$.  Then the generalization is
\begin{equation}\label{Eqn:BetaGamma}
\lab \beta_t \rab \leq \frac{\gamma_t}{2r_t^2} \quad\text{assuming that $r_t\neq 0$.}
\end{equation}
For completeness, we now take a moment to derive this inequality.  Because $\sigma_t$ is a projection, $\sigma_t^2=\sigma_t$.  Thus, in terms of a unit vector $(n_1,n_2,n_3)$ in the kernel of $\sigma_t$, the diffusion matrix associated to $Z_t$ is (and where we feel free to drop the $t$'s in the subscripts to simplify the notation)
\[\begin{bmatrix}
1-n_1^2 & -n_1n_2 & -n_1n_3 \\
-n_1n_2 & 1-n_2^2 & -n_2n_3 \\
-n_1n_3 & -n_2n_3 & 1-n_3^2
\end{bmatrix} . \]
The Hessian of $\log r$ (where $r\neq 0$) is given by
\begin{multline*}
\frac{\partial^2}{\partial x_1^2} \log r = \frac{1}{r^2}-\frac{2x_1^2}{r^4}, \quad
\frac{\partial^2}{\partial x_1\partial x_2} \log r= \frac{\partial^2}{\partial x_2\partial x_1} \log r
 = -\frac{2x_1x_2}{r^4}, \\
\frac{\partial^2}{\partial x_2^2} \log r = \frac{1}{r^2}-\frac{2x_2^2}{r^4}, \quad
\text{and $0$ for everything else.}
\end{multline*}
Then Ito's rule gives
\[
\beta_t = \frac{2x_1x_2n_1n_2}{r^4}+\frac{1}{2r^2}\lb \lp 1-\frac{2x_1^2}{r^2}\rp\lp1-n_1^2\rp
+ \lp 1-\frac{2x_2^2}{r^2}\rp\lp1-n_2^2\rp \rb .
\]
Similarly, we compute that $\gamma_t = 1-n_3^2 = n_1^2+n_2^2$.  Because the desired inequality is invariant under rotation around the $x_3$-axis, it is enough to prove it under the assumption that $x_1=r$ and $x_2=0$ at the current point.  Making this assumption, we get the simplification $\beta_t= (n_1^2-n_2^2)/2r^2$.  Comparing this with $\gamma_t/2r^2$ establishes the desired inequality.  In the case when $Z_t$ is Brownian motion on $M$, we have that $\beta_t = (\Lap_M/2)\log r$ and $\gamma_t = |\nabla_M x_3|^2$, where the right-hand sides are evaluated at $Z_t$, of course.  Thus the situation reduces to that of Lemma \ref{Lem:Basic}, interpreted probabilistically.

Next, our assumptions on $f$ (and the definition of the $\sigma_n$) imply that
\[\begin{split}
& \E\lb \left. \int_{\sigma_n}^{\sigma_{n+1}\wedge \zeta} \gamma_t \, dt 
\right| \{X_n=m\}\cap \sF^Z_{\sigma_n} \rb \\
& \quad\quad= \E\lb x_3^2\lp Z_{\sigma_{n+1}\wedge\zeta}\rp - x_3^2\lp Z_{\sigma_n}\rp \left| \{X_n=m\}\cap \sF^Z_{\sigma_n} \right. \rb \\
& \quad\quad\leq \lp f \lp e^{m+1}\rp \rp^2 .
\end{split}\]
Note that this bound is independent of everything in $\sF^Z_{\sigma_n}$ except $\{X_n=m\}$.  We use this estimate along with Equation \eqref{Eqn:BetaGamma} and Ito's rule to see that
\[
\E\lb \log r_{\sigma_{n+1}\wedge\zeta} | \{X_n=m\}\cap \sF^Z_{\sigma_n} \rb \leq n+\frac{\lp f \lp e^{m+1}\rp \rp^2}{2e^{2m-2}} ,
\]
which again is independent of everything in $\sF^Z_{\sigma_n}$ except $\{X_n=m\}$.  Recall that, given $X_n=m$, $X_{n+1}=\log r_{\sigma_{n+1}}$ is either $m+1$ or $m-1$, assuming $\sigma_{n+1}<\infty$ (that is, the process isn't stopped first), and if $\sigma_{n+1}=\infty$ then $\log r_{\zeta}\in[m-1,m+1]$.  Thus, we see that
\[
\Prob\lp X_{n+1}=m+1 | X_n=m\rp \leq \frac{1}{2}+\frac{\lp f \lp e^{m+1}\rp \rp^2}{4e^{2m-2}},
\]
independent of anything else in $\sF^Z_{\sigma_n}$.  In particular, we have found the desired uniform estimate for $\Prob\lp X_{n+1}=m+1 | X_n=m\rp$.

We wish to use this bound as the transition probability $p_m$ for our comparison chain.  This means that it must be less than one (recall that we want our chain to be irreducible).  Thus, we will now also assume that $L$ and $f(r)$ are such that
\begin{equation}\label{Eqn:fCondition}
\frac{\lp f \lp e^{m+1}\rp \rp^2}{4e^{2m-2}} <\frac{1}{2} \quad\text{for all $m\geq L+1$.}
\end{equation}
(Note that this is done mainly for convenience.)  We set
\[
p_m = \frac{1}{2}+\frac{\lp f \lp e^{m+1}\rp \rp^2}{4e^{2m-2}}\quad\text{for all $m\geq L+1$.}
\]

Next, we assume that the probability space on which $Z_t$ is defined is rich enough to support a countable sequence $U_1, U_2, \ldots$ of i.i.d.\ random variables, distributed uniformly on the interval $[0,1]$, and all independent of $Z_t$ (this can always be accomplished; for example, by taking the product of the original probability space with the unit interval with Lebesgue measure).  We now determine the process $Y_n$ as follows.  Let $Y_0=L+1$, and define $Y_n$ for $n\geq 1$ inductively as follows.  If $Y_{n-1}=m>X_{n-1}$ or if $X_{n-1}$ doesn't exist (equivalently, $\zeta<\sigma_{n-1}$), then
\[
Y_n = 
\left\{ \begin{array}{ll}
                m+1 & \text{if $U_n\in\lb 0,p_{m}\rb$,} \\
                m-1 & \text{if $U_n\in\lp p_{m},1\rb$.}
                \end{array}\right.  
\]
If  $Y_{n-1}=m=X_{n-1}$, then let
\[
\phi_m=\Prob\lb X_n=m+1 \left| \{X_{n-1}=m\}\cap \sF^Z_{\sigma_{n-1}}\right. \rb .
\]
Note that $\phi_m$ is bounded from above by $p_m$ (so $\phi_m\leq p_m$), and let
\[
Y_n = 
\left\{ \begin{array}{ll}
                m+1 & \text{if $X_n=m+1$,} \\
                m+1 & \text{if $X_n\neq m+1$ and $U_n\leq \lp p_m- \phi_m\rp/\lp1- \phi_m\rp$,} \\
                m-1 & \text{if $X_n\neq m+1$ and $U_n> \lp p_m- \phi_m\rp/\lp1- \phi_m\rp$.}
                \end{array}\right.  
\]
Here $X_n\neq m+1$ should be thought of as short-hand for $X_n=m-1$ or $X_n$ doesn't exist.  Finally, we stop $Y_n$ if it hits $L$.  Recall that we're also assuming that $Z_t$ is stopped if it hits $\{r=e^L\}$.

A simple computation shows that $\Prob\lb Y_n=m+1|Y_{n-1}=m\rb=p_m$, independent of anything else in $\sF^Z_{\sigma_{n-1}}$ and of $U_1,\ldots,U_{n-1}$.  It follows that $Y_n$ is a Markov chain with transition probabilities $p_m$.  In addition, by construction $Y_n$ dominates $X_n$ in the sense that $Y_n \geq X_n$ for all $n$ such that $\sigma_n\leq \zeta$, almost surely.  So $Y_n$ is the desired comparison chain.

\section{Parabolicity}

Because the conditions for minimal surfaces-with-boundary we are ultimately interested in, namely parabolicity and quadratic area growth, aren't affected by homothetically rescaling (as mentioned in the introduction), we will assume that $M$ is contained in the region $\{r\geq e^L\}$ for some non-negative integer $L$ (instead of the more general $\{r\geq c\}$ for positive $c$), and that any intersection of $M$ and $\{r=e^L\}$ is contained in the boundary of $M$.

Let
\[
f_1(r) = \frac{c r}{\sqrt{\log r \log\lp\log r\rp}} \quad\text{for some $c>0$},
\]
and let $L$ be chosen so that Equation~\eqref{Eqn:fCondition} is satisfied with $f=f_1$.  Note that a simple computation shows that such an $L$ always exists.  (We should perhaps use $L_1$ here, analogously to $f_1$, but that would make the subscripts too unwieldy.)

The point of working with rank 2 martingales in the last section is that we have the following more general version of parabolicity.

\begin{THM}\label{THM:MartPara}
Let $f_1(r)$ and $L$ be as in the preceding paragraph.  Suppose that $Z_t$ is a rank 2 martingale as described in Section \ref{Sect:Radial} (in particular, we continue to assume, for convenience, that $r(Z_0)=e^L$).  If $Z_t$ is contained in the region
\[
\lc r \geq e^L \text{ and } |x_3|\leq f_1(r) \rc \subset \Re^3 ,
\]
then $Z_t$ almost surely has finite lifetime.
\end{THM}
\emph{Proof:} We recall $X_n$, the discretization of $\log r(Z_t)$ introduced above.  Then with $f=f_1$, our earlier computations show that we have a comparison Markov chain $Y_n$, started at $Y_0=L+1$, with (after a bit of simplification)
\[
p_m= \frac{1}{2}+\frac{c^2 e^4}{4}\cdot\frac{1}{(m+1)\log(m+1)} .
\]
Then a little more algebra shows that
\[
p_m\leq \frac{m+1}{2m+1}\quad\text{if and only if} \quad
\frac{4}{c^2e^4}(m+1)\log(m+1)\geq 4m+2 .
\]
This last inequality holds for all sufficiently large $m$.  Thus $p_m\leq(m+1)/(2m+1)$ for large $m$, and it follows from the discussion in Section~\ref{Sect:MC} that $Y_n$ is parabolic, that is, $Y_n$ almost surely hits $L$.

Recall that $Y_n$ dominates $X_n$ until $\zeta$.  Then because $Z_t$ is contained in the region $\{r\geq e^L\}$ and we stop $Z_t$ when $\log r(Z_t)=L$, if not before, $Z_t$ is stopped no later than the time when $Y_n$ hits $L$, where this has the obvious meaning in terms of the stopping times $\sigma_n$.  In particular, $\zeta$ is almost surely finite.  $\Box$

Returning to minimal surfaces, recall that asserting that a surface-with-boundary is parabolic implies that the boundary is nonempty.  Again, we are interested in the asymptotic behavior of minimal surfaces, so the following theorem is formulated without much concern for the behavior of the minimal surface for ``small $r$.''

\begin{THM}\label{THM:Para}
Let $f_1(r)$ and $L$ be as above.  Suppose that $M$ is a stochastically complete minimal surface-with-boundary contained in the region
\[
\lc r \geq e^L \text{ and } |x_3|\leq f_1(r) \rc \subset \Re^3 ,
\]
with the assumption that that any intersection of $M$ and $\{r=e^L\}$ is contained in the boundary of $M$.  Then $M$ is parabolic.
\end{THM}
\emph{Proof:} We consider Brownian motion $B_t$ on $M$ started at some interior point with $r=e^{L+1}$, which we can assume exists after homothetically rescaling (though we might have to change $c$ and $L$ after rescaling) and using that $M$ is stochastically complete.  We let $\zeta$ be the first hitting time of the boundary and stop $B_t$ at $\zeta$.  Recall that, in order to show that $M$ is parabolic, it is enough to show that Brownian motion started at any interior point hits the boundary in finite time, that is, $\zeta$ is almost surely finite. Now $B_t=Z_t$ satisfies the hypotheses of Theorem \ref{THM:MartPara}, and we conclude that $B_t$ almost surely hits the boundary in finite time. $\Box$

The same arguments, along with more involved computations, could presumably be pushed to give a slightly weaker hypotheses, meaning replacing $f_1(r)$ with something larger by a sub-logarithmic factor.  However, these arguments, in their present form, break down if we try to extend them to the case of surfaces (or rank 2 martingales) contained in the exterior of a cone (that is, the case when $f_1(r)$ is replaced by $cr$ for some $c>0$).  The reason is that such a region is preserved under homothetic rescaling, and thus the $p_m$ (as defined above) don't decay to $1/2$ as $m\ra\infty$.  In another direction, one could estimate quantities associated with the parabolicity of $M$, for example, the distribution of the maximum of $r$ along Brownian paths.  However, we have not pursued these directions here because we are unaware of any applications of this kind of result.

The reason for framing the discussion in terms of rank 2 martingales is that, depending on your point of view, it provides an appealing ``explanation'' for parabolicity, and clarifies the hypotheses.  In particular, it is interesting that parabolicity doesn't rely on having a surface at all, or equivalently, on the plane field determining the evolution of the rank $2$ martingale being integrable.

\section{Area growth}

While it is common to study parabolicity using Brownian motion, it is perhaps more interesting that it also makes a fairly natural tool, at least in this context, to study area growth.  The generalization of quadratic area growth to rank 2 martingales would be a quadratic estimate on occupation times, as the proof the following theorem makes clear.  However, working with rank 2 martingales doesn't seem worth the added complexity in the present context.  Thus, for the remainder of the section we restrict our attention to the case when our rank 2 martingale is Brownian motion on a stochastically complete minimal surface-with-boundary, and adopt notation accordingly.

Let
\[
f_2(r) = \frac{c r}{\sqrt{\log r} \log\lp\log r\rp} \quad\text{for some $c>0$},
\]
and let $L$ be chosen so that Equation~\eqref{Eqn:fCondition} is satisfied with $f=f_2$.  Again, a simple computation shows that such an $L$ always exists.  Note that for a minimal surface satisfying the assumptions of the following theorem, quadratic area growth is equivalent to quadratic area growth with respect to $r$, as discussed in the introduction.

\begin{THM}\label{THM:Area}
Let $f_2(r)$ and $L$ be as in the preceding paragraph.  Suppose that $M$ is a properly immersed minimal surface with compact, non-empty boundary, contained in the region
\[
\lc   r\geq e^L \text{ and } |x_3|\leq f_2(r)  \rc \subset \Re^3 
\]
with the assumption that any intersection of $M$ and $\{r=e^L\}$ is contained in the boundary of $M$.  Then $M$ has quadratic area growth.
\end{THM}

\emph{Proof:} 
Removing a compact subset of $M$ doesn't affect the result, nor does homothetic rescaling (though we might have to change $c$ and $L$ after rescaling).  Thus, without loss of generality, we can assume that the boundary of $M$ is equal to $M\cap \{r =e^L\}$ and that $M\cap \{r= e^{L+1}\}$ is a smooth, compact curve in $M$, which we denote $\gamma$.  Then $\gamma$ divides $M$ into two connected components (which intersect only at $\gamma$), the unbounded piece $M_u=M\cap \{r\geq e^{L+1}\}$ with boundary $\gamma$, and the compact piece $M_c=M\cap \{e^L\leq r\leq e^{L+1}\}$ with two boundary components, $\gamma$ and $\partial M$ (the original boundary of $M$).  Also, since $\partial M= M\cap \{r =e^L\}$, we see that Brownian motion is only stopped when it hits the level $\{r=e^L\}$.

Note that for large $r$, $f_2(r)\leq f_1(r)$, regardless of what constants $c$ are used in the two functions.  So applying Theorem~\ref{THM:Para} and using that $M$ is properly immersed, we see that $M$ is parabolic.  It follows that $M_u$ is also parabolic.

We consider the function $h$ on $M$, determined as follows.  We let $h$ be harmonic on $M_c$ with boundary values 1 on $\gamma$ and 0 on $\partial M$.  Because $M_c$ is compact, these boundary values uniquely determine a function.  On $M_u$, we let $h$ be the unique bounded harmonic function with boundary values 1 on $\gamma$ (which is well-defined because $M_u$ is parabolic); thus $h$ is identically 1 on $M_u$.  Finally, we let $h$ be 1 on $\gamma$.  Thus, $h$ is a continuous, bounded function on $M$ and is harmonic on $M\setminus\gamma$.

Observe that $\frac{1}{2}\Lap h$ is measure-valued, supported on $\gamma$, and has finite mass; let that mass be $1/\alpha$.  Then $\frac{1}{2}\Lap(\alpha h)$ is a probability measure, which we denote $\mu$, supported on $\gamma$.  It follows that $\alpha h$ is the Green's function associated to $\mu$, and thus $\alpha h$ is also the density of the expected occupation time of Brownian motion on $M$ started from $\mu$ and stopped at $\partial M$.  (See \cite{Grigoryan} for background on the relationship between Brownian motion and potential theory on manifolds.)  Since we know that $M_c$ has finite area and $\alpha h$ is bounded there and that $\alpha h$ is identically equal to $\alpha$ on $M_u$, it follows that the expected occupation time of $M\cap\{ r\leq \rho\}$ is asymptotic to $\alpha$ times the area of $M\cap\{ r\leq \rho\}$ (as $\rho\ra\infty$).  This is the connection between Brownian motion and area growth that we need.  In particular, suppose that, for Brownian motion started at some point on $\gamma$, we show that the expected occupation time of $M\cap\{ r\leq \rho\}$ is bounded from above by $C\rho^2$ for some positive constant $C$. If this holds with $C$ independent of which point of $\gamma$ our Brownian motion starts from, then the estimate will hold for Brownian motion started at $\gamma$ and we will have shown that $M$ has quadratic area growth.

To do this, let $\rho=e^k$ for some integer $k>L+1$.  We start be observing that Brownian motion started at a point with $r=e^{L+1}$ (that is, a point on $\gamma$) almost surely hits either $r=e^L$ (in which case the process is stopped) or $r=e^{k+1}$ (that is, $e\rho$) in finite time.  Next, we need to estimate how long this takes.  Let $\theta_j$ be first hitting time of $r=e^{j}$ for any integer $j$; then 
\[
\E\lb r^2_{\theta_L\wedge\theta_{k+1}} |r_0=e^{L+1} \rb =e^{2(k+1)}\Prob\lp\theta_{k+1}<\theta_L \rp
+ e^{2L}\Prob\lp\theta_L<\theta_{k+1} \rp
\]
(independent of where on $\{r=e^{L+1}\}$ the Brownian motion starts).  On the other hand, from $1\leq (\Lap_M /2)r^2\leq 2$ and Ito's formula, we see that
\[
\E\lb r^2_{\theta_L\wedge\theta_{k+1}} |r_0=e^{L+1} \rb \geq e^{2(L+1)} + \E\lb \theta_L\wedge\theta_{k+1}\rb
\]
(again independent of where the Brownian motion starts).  It follows that
\[
\E\lb\theta_L\wedge\theta_{k+1} | r_0=e^{L+1}\rb \leq e^{2(k+1)}-e^{2(L+1)} .
\]
This is the desired estimate on how long it takes for Brownian motion to leave $r\in(e^L,e^{k+1})$; note that it is quadratic in the upper endpoint $e^{k+1}=e\rho$ and independent of where on $\gamma=\{r=e^{L+1}\}$ the Brownian motion starts.

Continuing, we see that if $\theta_L<\theta_{k+1}$ then the process is stopped at $\theta_L$ and the only contribution to the occupation time of $M\cap\{ r\leq \rho\}$ is $\theta_L$.  (Because we're stopping at $r=e^{k+1}$ rather than $e^k$, this is just an upper bound on the contribution; the reason for this set-up will be clear below.)  If $\theta_{k+1}<\theta_L$ then the process continues, and we must also estimate the contribution to the occupation time after $\theta_{k+1}$.  Intuitively, we do this as follows.  Once the process hits $r=e^{k+1}$, we stop counting it toward the occupation time, until the next time it hits $r=e^k$ (which almost surely happens).  Then we start counting the time again, until the next time it hits either $e^L$ or $e^{k+1}$, at which point either the process is stopped or we perform another round of the same procedure.  Note that, almost surely, the process eventually hits $r=e^L$ and is stopped.  This procedure counts all of the occupation time of $M\cap\{ r\leq \rho\}$ and part of the occupation time of $M\cap\{ \rho < r\leq e^{k+1}\}$, but it certainly provides an upper bound on the occupation time of $M\cap\{ r\leq \rho\}$.  Further, for the purpose of establishing quadratic area growth, this is no real loss, since $e^{L+1}=e\cdot\rho$.

We now formalize the above procedure.  Let $\xi_0 = \theta_{k+1}$ be the first hitting time of $r=e^{k+1}$, and let $\tilde{\xi}_0$ be the first hitting time of $r=e^{k}$ after $\xi_1$.  Then we inductively define $\xi_n$ to be the first hitting time of $r=e^{k+1}$ after $\tilde{\xi}_{n-1}$ and $\tilde{\xi}_n$ to be the first hitting time of $r=e^{k}$ after $\xi_n$.  By parabolicity, $\tilde{\xi}_n$ is almost surely finite if $\xi_n$ is.  Let $U_k=\sup\{ n\geq 0: \xi_n<\infty\}$ be the number of upcrossings from $r=e^k$ to $e^{k+1}$; again by parabolicity it is almost surely finite.  An argument completely analogous to the above shows that, starting from $\tilde{\xi}_{n-1}$, the expected amount of time spent until the process leaves $r\in(e^L,e^{k+1})$ satisfies the estimate
\[
\E\lb\lp\theta_L\wedge\xi_n\rp - \tilde{\xi}_{n-1}\right|\left. \tilde{\xi}_{n-1}<\infty \rb \leq e^{2(k+1)}-e^{2k} .
\]
(We note that this estimate is uniform over where on $\{r=e^k\}$ the process is at $\tilde{\xi}_{n-1}$.) Thus, the contribution of each upcrossing to the expected occupation time is bounded by the expression on the right.

Summarizing our progress so far, we have that the expected occupation time of $M\cap\{ r\leq \rho\}$, where $\rho=e^k$, for Brownian motion started at any point of $\gamma$ is bounded from above by
\[
e^{2(k+1)}-e^{2(L+1)} + \Prob\lp \theta_{k+1}<\theta_L\rp\lb e^{2(k+1)}-e^{2k}\rb
\E\lb U_k |  \theta_{k+1}<\theta_L \rb .
\]
Observe that $e^{2(k+1)}-e^{2(L+1)}$ and $e^{2(k+1)}-e^{2k}$ are both quadratic in $\rho=e^k$.  Thus, in order to show that this expression grows quadratically (in $\rho$), it is enough to show that
\[
\Prob\lp \theta_{k+1}<\theta_L\rp \E\lb U_k |  \theta_{k+1}<\theta_L \rb
\]
is bounded as $k\ra\infty$.  Again, we recall $X_n$, the discretization of $\log r_t$ introduced above.  We see that $\theta_{k+1}$ corresponds to the first time $X_n$ hits $k+1$ (where this correspondence is understood in the natural way in terms of the stopping times $\sigma_n$) and $\theta_L$ corresponds to the first time $X_n$ hits $L$ (and is stopped).  Also, note that $\zeta=\theta_L$, so $X_n$ is a nearest-neighbor random walk on $\{L,L+1,\ldots\}$ stopped when, and only when, it hits $L$.  Thus we will have no need to refer explicitly to $\zeta$, and the discussion in Section \ref{Sect:Radial} simplifies accordingly in this case.  (This simplification is made possible because here, unlike in Theorem \ref{THM:Para}, we are assuming that $M$ is properly immersed with compact boundary, as discussed at the beginning of the proof.)

We see that
$\Prob\lp \theta_{k+1}<\theta_L\rp$ can be re-expressed as  the probability that $X_n$ hits $k+1$ before $L$.  Further, letting $f=f_2$, our earlier computations show that we have a comparison Markov chain $Y_n$ with
\[
p_m= \frac{1}{2}+\frac{c^2e^4}{4}\cdot\frac{1}{(m+1)\lp \log(m+1)\rp^2} . 
\]
Thus $\Prob\lp \theta_{k+1}<\theta_L\rp$ is bounded from above by the probability that $Y_n$, started from $L+1$, hits $k+1$ before $L$, independent of where on $\gamma$ the underlying Brownian motion in $M$ starts.  Further, this probability was given in Section \ref{Sect:MC} as the reciprocal of $A_k$, and recalling that gives
\[
\Prob\lp \theta_{k+1}<\theta_L\rp \leq 
\lp 1+\frac{q_{L+1}}{p_{L+1}} +\frac{q_{L+1}q_{L+2}}{p_{L+1}p_{L+2}} +\cdots +\frac{q_{L+1}\cdots q_k}{p_{L+1}\cdots p_k} \rp^{-1},
\]
where $p_m$ is as above and $q_m=1-p_m$, as usual.  Also, observe that both $X_n$ and $Y_n$ are parabolic, which follows from the proof of Theorem \ref{THM:Para}.

Next, we wish to estimate $\E\lb U_k |  \theta_{k+1}<\theta_L \rb$ in terms of $Y_n$ as well, and note that such an estimate in terms of $Y_n$ will hold for $X_n$ independent of where on $\gamma$ the underlying Brownian motion starts.  The point is that, even though upcrossings of $X_n$ might occur at different times than upcrossings of $Y_n$, parabolicity and the Markov property (of $Y_n$) imply that $\E\lb U_k |  \theta_{k+1}<\theta_L \rb$ is no larger than $\E\lp u_k| Y_0=k \rp$.  This is more or less obvious, but we provide a bit more detail for the sake of clarity.  The main idea is that we can feel free to pause either $X_n$ or $Y_n$ in order to allow the other process to ``catch up'' without changing the distribution of the total number of upcrossings of either process.  First, because we're conditioning on $\theta_{k+1}<\theta_L$, we know that $X_n$ will hit $k+1$ and then by parabolicity will come back to $k$.  We wait until the happens to start $Y_n$ (which we start from $k$).  Then $Y$ dominates $X$, so if $X$ has an upcrossing (that is, if $X$ hits $k+1$ again before hitting $L$), so does $Y$.  Supposing that $X$ does have an upcrossing, it might come back down to $k$ before $Y$ does.  However, once $X$ comes back down to $k$, we can ``pause'' it there and wait for $Y$ to also come back down to $k$, which it will by parabolicity.  With both processes at $k$, we simply repeat the above procedure.  This shows that every upcrossing of $X$ corresponds to an upcrossing of $Y$, and thus the desired inequality for the expected numbers of upcrossings holds.

Now using our computation of $\E\lp u_k| Y_0=k \rp$ from Section \ref{Sect:MC}, we see that
\begin{multline*}
\Prob\lp \theta_{k+1}<\theta_L\rp \E\lb U_k |  \theta_{k+1}<\theta_L \rb \leq
\lb 1+\frac{p_k}{q_k}+\frac{p_k p_{k-1}}{q_k q_{k-1}}+\cdots+\frac{p_k\cdots p_{L+1}}{q_k\cdots q_{L+1}}\rb \\
\times \lb 1+\frac{q_{L+1}}{p_{L+1}} +\frac{q_{L+1}q_{L+2}}{p_{L+1}p_{L+2}} +\cdots +\frac{q_{L+1}\cdots q_k}{p_{L+1}\cdots p_k}\rb^{-1} .
\end{multline*}
Recall that, in order to show that the expected occupation time of $M\cap\{ r\leq \rho\}$, where $\rho=e^k$, grows quadratically in $\rho$, it is enough to show that the above quantity stays bounded as $k\ra\infty$.  The right-hand side is a fraction (though not displayed as such), and because $p_m>q_m$ the terms in the numerator are increasing (from left to right) while the terms in the denominator are decreasing.  Since both numerator and denominator have the same number of terms (for a given $k$), we see that in order to show that this fraction is bounded as $k\ra\infty$, it is sufficient to show that the ratio of the right-most terms, that is,
\[
\lp \frac{p_k\cdots p_{L+1}}{q_k\cdots q_{L+1}}\rp\Biggl/
\lp \frac{q_{L+1}\cdots q_k}{p_{L+1}\cdots p_k}\rp = \lp \frac{p_k\cdots p_{L+1}}{q_k\cdots q_{L+1}}\rp^2
\]
is bounded as $k\ra\infty$.  Taking the logarithm of the right-hand side, we see that it is enough to show that $\sum_{m=L+1}^k\log(p_m/q_m)$ is bounded as $k\ra\infty$ (recall that all terms in the sum are positive since $p_m>q_m$).  To do this, we will show that
\[\tag{$*$}
\log\lp\frac{p_m}{q_m}\rp \leq \frac{\tilde{c}}{(m+1)\lp \log(m+1)\rp^2}
\]
for some $\tilde{c}>0$, for large $m$ (so that $\sum_{i=L+1}^k\log(p_m/q_m)$ is bounded by the usual comparison test for sums).

We have (using the above formula for $p_m$) that
\[\begin{split}
\log\lp\frac{p_m}{q_m}\rp = &\log\lp \frac{1}{2}+\frac{c^2e^4}{4}\cdot\frac{1}{(m+1)\lp \log(m+1)\rp^2}\rp \\
&- \log\lp \frac{1}{2}-\frac{c^2e^4}{4}\cdot\frac{1}{(m+1)\lp \log(m+1)\rp^2}\rp .
\end{split}\]
For large $m$, both $p_m$ and $q_m$ become arbitrarily close to $1/2$.  Thus we can use an estimate based on the first-order Taylor expansion of the logarithm, namely that $|\log(x+1/2)-\log(1/2)|\leq 3|x|$ for $x$ sufficiently close to 0.  This implies that, for large $m$,
\[\begin{split}
\log\lp\frac{p_m}{q_m}\rp &\leq \log\lp \frac{1}{2}\rp +\frac{3c^2e^4}{4}\cdot\frac{1}{(m+1)\lp \log(m+1)\rp^2} \\
&\quad\quad- \log\lp \frac{1}{2}\rp +\frac{3c^2e^4}{4}\cdot\frac{1}{(m+1)\lp \log(m+1)\rp^2} \\
& = \frac{3c^2e^4}{2}\cdot\frac{1}{(m+1)\lp \log(m+1)\rp^2} .
\end{split}\]
This establishes $(*)$.  Combining all of the above, we have proven the desired estimates on Brownian motion uniformly with respect to where on $\gamma$ the Brownian motion starts, and we can now conclude that there exists a positive constant $C$ such that the area of $M\cap\{ r\leq \rho\}$, when $\rho=e^k$, is bounded from above by $C\rho^2$ for sufficiently large $\rho$.

Finally, to extend this result to all $\rho$ (and not just those that are an integer power of $e$, as we have been assuming to make use of our discretization $X_n$ of the radial process), suppose that $\rho\in(e^{k-1},e^k]$.  Then since the area is monotone in $\rho$, we have
\[
\Area\lp M\cap\{ r\leq \rho\}\rp \leq Ce^{2k} =Ce^2e^{2(k-1)}<\lp Ce^2\rp \rho^2 .
\]
So increasing $C$ by a factor of $e^2$ gives the estimate for general $\rho$, completing the proof.  $\Box$

Note that, having proven the above result, a standard argument can be used to extend it somewhat.  In particular, the monotonicity formula for area (see \cite{Simon}) shows that
\[
\frac{\Area\lp M\cap\{ r\leq \rho\}\rp}{\pi\rho^2}
\]
converges to a finite limit as $\rho\ra\infty$, which we denote $n(M)$.  Then because $M$ is contained in
\[
\lc  r\geq e^L \text{ and } |x_3|\leq f_2(r)  \rc ,
\]
geometric measure theory allows us to conclude that the homothetic shrinkings  of $M$ converge to a locally finite minimal varifold supported on the $(x_1,x_2)$-plane.  This limit is thus an integer multiple of the plane, and we conclude that $n(M)$ is a positive integer.  (See either the end of Section 2 in \cite{CKMR} or Section 6.4 of \cite{MPSurvey} for more details of the argument.)  When $M$ is an appropriate representative of an end of a minimal surface, $n(M)$ is generally called the multiplicity of the end.

Recall that if the area of a surface grows quadratically with respect to the intrinsic distance, then the surface is parabolic.  In light of this, it is not surprising that these methods give parabolicity under a strictly weaker condition than is needed for quadratic area growth (with respect to $r$, or equivalently the exterior distance from the origin).  Also, as was the case for Theorem~\ref{THM:Para}, these computations could likely be extended to allow replacing $f_2(r)$ by something slightly more lenient, but they cannot be simply extended to allow $f_2(r)$ to be replaced by $f_1(r)$, much less to allow $M$ to be contained in the complement of a cone.  Further, in various cases one could give explicit estimates for the constants appearing in the proof in order to derive an estimate for $C$, and thus also for $n(M)$.  In particular, $\alpha$ depends only on the behavior of a compact subset of $M$, such as $M\cap \{r\leq e^{L+2}\}$.  Moreover, $\Prob\lp \theta_{k+1}<\theta_L\rp \E\lb U_k |  \theta_{k+1}<\theta_L \rb$ can be estimated in terms of the $p_m$, which depend only on $f_2$.  Replacing $f_2$ with a smaller $f$ would decrease the $p_m$ and thus also the upper bound on $n(M)$.  However, as was the case previously, we have not pursued this because we are unaware of any interest in such an estimate.

\bibliographystyle{amsplain}

\def\cprime{$'$}
\providecommand{\bysame}{\leavevmode\hbox to3em{\hrulefill}\thinspace}
\providecommand{\MR}{\relax\ifhmode\unskip\space\fi MR }
\providecommand{\MRhref}[2]{%
  \href{http://www.ams.org/mathscinet-getitem?mr=#1}{#2}
}
\providecommand{\href}[2]{#2}

\end{document}